\documentclass[12pt]{amsart}

\usepackage{amssymb}
\usepackage[centertags]{amsmath}
\usepackage{amscd, amsthm}
\usepackage{stmaryrd}
\usepackage{enumitem}

\setlist[enumerate,1]{label={(\arabic*)}} 

\usepackage{xspace}

\usepackage{hyperref}

\newtheorem{prop}{Proposition}[section]

\newtheorem{thm}[prop]{Theorem}

\newtheorem{defini}[prop]{Definition}
\newtheorem{osserva}[prop]{Remark}

\newtheorem*{theorem*}{Theorem}

\theoremstyle{remark}

\def\be{\begin{equation}}
\def\ee{\end{equation}}

\def\dotminus{\mathbin{\ooalign{\hss\raise1ex\hbox{.}\hss\cr
  \mathsurround=0pt$-$}}}

\begin{document}
\title{Products of pseudofinite structures}

\author[P. D'Aquino]{Paola D'Aquino}
\address{Dipartimento di  Matematica e Fisica\\ Universit\`a della Campania L.Vanvitelli\\
Viale Lincoln 5, 81100 Caserta, Italy}
\email{paola.daquino@unicampania.it}

\author[A.J. Macintyre]{Angus J. Macintyre}
\address{School of Mathematics\\ University of Edinburgh\\ EH9 3FD Edinburgh, UK} 
\email{a.macintyre@qmul.ac.uk}
\subjclass{03C20, 03C52}
\keywords{Pseudofinite structure, product of first order structures, Boolean algebras.} 
\thanks{First author partially supported by PRIN 2022 ‘‘Models, Sets, and Classifications"}
\maketitle

\begin{abstract}
We prove that any product of a family of  pseudofinite structures is pseudofinite. The main tools are the fundamental results on products of first order structures in \cite{FV}.

\end{abstract}

\section{Introduction}


Fix a first order language $\mathcal L$. 

\begin{defini}
\label{pseudo}
An   $\mathcal L$-structure $\mathcal M$ is pseudofinite if  every $\mathcal L$-sentence true in $\mathcal M$  is true in some finite $\mathcal L$-structure. 
\end{defini}

We recall the following characterizations of pseudofinite structures. 

\begin{prop}
Let $\mathcal L$ be a 1st order language, and $\mathcal M$ an $\mathcal L$-structure. The following are equivalent
\begin{enumerate}
\item 
$\mathcal M$ is pseudofinite
\item
$\mathcal M$ is elementarily equivalent to an ultraproduct of finite $\mathcal L$-structures
\item
$\mathcal M\models Fin_{\mathcal L}$, where $ Fin_{\mathcal L}=\{\varphi : \mathcal N\models \varphi \mbox{ for all finite } {\mathcal L}-\mbox{structures } \mathcal N \}$.
\end{enumerate}
\end{prop}

 Note that a finite structure is pseudofinite. However, Ax in a great paper  \cite{ax4}  (preceding the study of pseudofinite structures) defined for fields, a different notion of pseudofinite, that agrees with the general notion for infinite fields, while finite fields are not pseudofinite  in the sense of Ax. For fields, Ax's notion is directly linked to Weil's  Riemann Hypothesis for curves, and it is this that makes it revolutionary. 
 
 In this paper, we will adopt the convention of using the notion of pseudofinite which includes finite structures. This excludes Ax's notion, but in fact everything we do here works for Ax's notion, whose models are exactly the infinite fields which are pseudofinite in the general sense. 
 
 Examples of pseudofinite structures are $(\mathbb Q,+)$,  and random graph, and examples of non pseudofinite structures are $(\mathbb Q,<)$, $(\mathbb Z,+)$, and any algebraically closed field. 
It is an easy exercise to show that pseudofiniteness is preserved under interpretability, reducts  and expansions of a language. It is trivial that pseudofiniteness is preserved under ultraproducts, but in general it is not preserved under reduced product as we will show in the final remark. 

For elementary properties of pseudofinite structures we refer to the well written notes by D. Garcia \cite{Garcia}.

Pseudofinite rings play a crucial role in the model theoretic analysis of the rings $\mathcal M/n\mathcal M$ where $\mathcal M$ is a non standard model of Peano Arithmetic, that the authors conduct in \cite{resrings} and in \cite{CommUnital}. 

\section{Main result}
Our goal is to show that product of pseudofinite structures is pseudofinite. If $(\mathcal M_i)_{i\in I}$ are $\mathcal L$-structures then $\prod_{i\in I}\mathcal M_i$ denotes the product of the  $\mathcal M_i$'s for $i\in I$, and  the induced $\mathcal L$-structure  is defined coordinatewise.  The main result in this short paper  is the following theorem.

\begin{thm}
\label{productpseudo}
Let $I$  be an index set (either finite or infinite).  
If   $(\mathcal M_i)_{i\in I}$ are pseudofinite $\mathcal L$-structures then $\prod_{i\in I}\mathcal M_i$ is pseudofinite. 
\end{thm}

For the proof we will use a classical result due to Feferman and Vaught (Theorem 3.1 in \cite{FV}). This was crucial also for establishing the connection between the residue rings of a model of $PA$ given by composite  elements and those given by prime powers, as proved by the authors in \cite{CommUnital}. This will be further developed  in the next and last section of this paper. 

Here we recall only the notions needed in the part of Theorem 3.1 in \cite{FV} that we will use, paying careful attention to uniformities. The language $\mathcal L$ is fixed while the index set $I$ is allowed to vary, as is the assignment  $i\mapsto \mathcal M_i$ for $i\in I$, and  $\mathcal L$-structures $\mathcal M_i$. The Boolean algebra $\mathcal P(I)$, the power set of the index set $I$, is an atomic Boolean algebra whose atoms are the singletons $\{i\}$ for all $i\in I$. Tarski showed that there is a uniform quantifier elimination for the theory of infinite atomic Boolean algebras  in the   language of Boolean algebras $\mathcal L_{B}= \{ \pmb{\wedge}, \pmb{\vee}, \pmb{\neg}, 0,1\} $ expanded by  infinitely many  unary predicates $C_j(x)$  saying that there are at least $j$ atoms below $x$.    Hence,  $$\mathcal P(I)\models C_j(w) \mbox{ if and only if  } w \subseteq I \mbox{ contains    at least } j \mbox{ elements  }.$$

It will be convenient for us to work in the language of Boolean rings where the binary ring operations $\oplus $ and $\odot $ are quantifier-free definable from $\pmb{\wedge}, \pmb{\vee}, \pmb{\neg}$ as follows 

\medskip

$
\begin{array}{ccc}
x\oplus y &   =    & (x\pmb{\wedge} \pmb{\neg}y)\pmb{\vee} (\pmb{\neg}  x\pmb{\wedge} y)  \\
x\odot y   &    =   &  x \pmb{\wedge} y   \\
0    &    =    &    0  \\
1    &    =    &    1    
\end{array}
$

\medskip

We will use the same notation for the two languages of Boolean algebras and of Boolean rings, i.e. $\mathcal L_{B}$. 

For any $I$ and $J$ infinite sets the Boolean algebras $\mathcal P(I)$ and $\mathcal P(J)$ are elementarily equivalent in the language $\mathcal L_{B}^C= \mathcal L_{B}  \cup \{ C_j(x): j\in \mathbb N, j>0\}$. The finite Boolean algebras are exactly the $\mathcal P(I_n)$ where $I_n=\{ 0, \ldots , n-1\}$, of cardinality $2^n$, and any non principal ultraproduct of them all is elementarily equivalent to the $\mathcal P(J)$ for some infinite set $J$.  It follows that for any $\mathcal L_{B}^C$-formula $\delta$ there is an integer $k=k(\delta)$ such that either

-  for all finite $I_n$ with $n\geq k(\delta) $, $\mathcal P(I_n)\models \delta$, or 

- for all finite $I_n$ with $n\geq k(\delta) $, $\mathcal P(I_n)\models \neg \delta$.

\medskip

By Tarski's elimination of quantifiers, for any $\mathcal L_{B}$-formula $\psi(\bar{x})$ there is a quantifier-free formula ${\psi^*}(\overline{x})$ in $\mathcal L_{B}^C$, such that if $I$ is infinite then $\mathcal P(I)\models \forall \bar{x}(\psi(\bar{x}) \leftrightarrow {\psi^*}(\bar{x}) )$.  Notice that ${\psi^*}_n$ contains at most the predicates $C_0(x), \ldots , C_n(x)$.

Moreover, there is $k=k(\psi)$ such that  for $n\geq k(\psi)$, in the finite Boolean algebra with $2^n$ elements the following equivalence holds,  $\mathcal P(I_n) \models \forall \bar{x}(\psi(\bar{x}) \leftrightarrow {\psi^*}(\bar{x}) )$.

For $n<k(\psi)$, the Boolean algebra  $\mathcal P(I_n)$   is ultrahomogeneous  in the language $\mathcal L_{B}^C$, and so it has elimination of quantifiers in $\mathcal L_{B}$  (see \cite{hodges}). There is a quantifier-free $\mathcal L_{B}$-formula ${\psi^*}_n(\bar{x})$ such that 
$ \mathcal P(I_n)\models \forall \bar{x}(\psi(\bar{x}) \leftrightarrow {\psi^*}_n(\bar{x}) ).$  

Note that the case of $I$ finite in the proof of Theorem \ref{productpseudo} is easy, and does not need to appeal to the elimination of quantifiers for atomic Booean algebras. 

\

We now combine  the elimination of quantifiers for atomic Boolean algebras with  Theorem 3.1 of \cite{FV}, and set up the combinatorics needed to prove Theorem \ref{productpseudo}.  In \cite{FV} the following notions of acceptable sequence and partition sequence were introduced. We recall these notions only for sentence since in what follows we only need them for sentences.

Let $(\mathcal M_i)_{i\in I}$ be $\mathcal L$-structures.  We expand $\mathcal L$ by another sort $\mathcal B$, the Boolean sort, and we will work with the 2-sorted language $\mathcal L_{\mathcal B}=\mathcal L \cup \mathcal L_{B}^C$, where $\mathcal L_{B}^C= \mathcal L_{B}  \cup \{ C_j(x): j\in \mathbb N, j>0\}$. The Boolean part of $\mathcal L_{\mathcal B}$ will be interpreted in the atomic  Boolean algebra $\mathcal P(I)$.

\begin{enumerate} 
\item An  {\it acceptable  sequence} is a $\zeta =\langle  \Psi, \theta_0, \ldots, \theta_m\rangle$ where  $\Psi $ is a $\mathcal L_{B}$-formula with at most $x_0, \ldots ,x_m$ free variables, and $\theta_0, \ldots, \theta_m$ are $\mathcal L$-sentences. 


\item  A sequence  $\langle \theta_0, \ldots, \theta_m\rangle$ of $\mathcal L$-sentences is called a {\it partition sequence } if 

\begin{enumerate}
\item $\theta_1\vee \ldots \vee \theta_m$ is true in all $\mathcal L$-structures, and

\item if $k\not=h$ then $\theta_k\wedge  \theta_h$ is false in all $\mathcal L$-structures.
\end{enumerate}

\item If $(\mathcal M_i)_{i\in I}$ is a family of $\mathcal L$-structures then for any sentence $\theta $ in $\mathcal L$  we set $$\llbracket \theta \rrbracket_I= \{i\in I: \mathcal  M_i\models \theta \}\in \mathcal P(I).$$

\end{enumerate}
For any fixed $\mathcal L$-formula $\theta $, $\llbracket \theta \rrbracket_I$  depends on the set $I$ and on the structures $\mathcal M_i$'s, only. The $\llbracket \cdot \rrbracket_I$ respects the Boolean algebra operations as specified in Section 2.3 of \cite{CommUnital}. Notice also that if $J, J_1,J_2\subseteq I$ then $\llbracket \theta \rrbracket_J\in \mathcal P(J)$, and $\llbracket \theta \rrbracket_{J_1} \cap \llbracket \theta \rrbracket_{J_2} = \llbracket \theta \rrbracket_{J_1\cap J_2}\in \mathcal P( J_1\cap J_2)$. In the sequel if no confusion can arise we will drop the subscript $I$.

If the sequence,  $\langle \theta_0, \ldots, \theta_m\rangle$ of $\mathcal L$-sentences is a partition then 

\begin{itemize}
\item 
 $\llbracket \theta_0 \rrbracket_I\oplus \ldots \oplus \llbracket \theta_m \rrbracket_I=I $

\item
$\llbracket \theta_h \rrbracket_I \odot \llbracket \theta_k \rrbracket_I=\emptyset \mbox{ for  }h,k\in \{0, \ldots, m\} \mbox{ and }h\not= k.$
\end{itemize}

We do not exclude that $\llbracket \theta_h \rrbracket_I =\emptyset$  for some $h\in \{0, \ldots , m\}$.

\smallskip

The following result is the content of Theorem 3.1 of \cite{FV}, relative only to sentences in the starting language $\mathcal L$. 

\begin{thm}
\label{FVsentences}
For each $\mathcal L$-sentence $\phi$ there is an acceptable sequence $$\zeta =\langle  \Psi, \theta_0, \ldots, \theta_m\rangle,$$ where $\Psi(y_0,\ldots ,y_m)$ is    $\mathcal L_{B}$-formula,  $\langle  \theta_0, \ldots, \theta_m\rangle$ is a partition sequence of  $\mathcal L$-sentences and for any index set $I$ and any family $(\mathcal M_i)_{i\in I}$ of $\mathcal L$-structures the following holds:
$$ \prod_{i\in I}\mathcal M_i \models \phi \mbox{ \hspace{.2in} iff \hspace{.2in} } \mathcal P(I) \models \Psi (\llbracket \theta_0 \rrbracket_I, \ldots ,\llbracket \theta_m \rrbracket_I).$$
\end{thm}

From the above discussion $\Psi(y_0,\ldots , y_m)$ is equivalent over $\mathcal P(I)$ for sufficiently {\it large} $I$ to a quantifier-free $\mathcal L_{B}^C$-formula $\Psi^*(y_0,\ldots , y_m)$ with at most $y_0,\ldots , y_m$ as free variables, or to an $\mathcal L_{B}$-formula $\Psi_n^*(y_0,\ldots , y_m)$   over $\mathcal P(I_n)$ for {\it small} $n$. Without loss of generality we can assume that the formula $\Psi^*(y_0,\ldots , y_m)$ is a conjunction of finitely many $C_q(f_h)$ and $\neg C_r(g_t)$,  where $f_h$ and $g_t$ are linear polynomials in the variables $y_0,\ldots , y_m$ (recall that in a Boolean ring $y^2=y$) over $\mathbb F_2$, and constant term either $0$ or $1$. We have further simplifications in the polynomials  $f_h$'s and $g_t$'s  when we evaluate $\Psi^*(y_0,\ldots , y_m)$ at $(\llbracket \theta_0 \rrbracket_I, \ldots ,\llbracket \theta_m \rrbracket_I)$. Indeed, from  $\langle \theta_0, \ldots, \theta_m\rangle$ being a partition it follows that for any $i\in I$ there is a unique $j\in \{0, \ldots , m\}$ such that $\mathcal M_i\models \theta_j$. So, $\llbracket \theta_j \rrbracket_I \odot \llbracket \theta_k \rrbracket_I=0$ for $j,k\in \{ 0, \ldots , m\}$ and $j\not=k$. This implies that the polynomials $f_h$ and $g_t$ evaluated at $\llbracket \theta_0 \rrbracket_I, \ldots ,\llbracket \theta_m \rrbracket_I$ are elements of $\mathcal P(I)$ of the form either 

\begin{enumerate}
\item
$\epsilon_0\llbracket \theta_0 \rrbracket_I\oplus \ldots \oplus \epsilon_m\llbracket \theta_m \rrbracket_I$ or

\item
$\epsilon_0\llbracket \theta_0 \rrbracket_I\oplus \ldots \oplus \epsilon_m\llbracket \theta_m \rrbracket_I\oplus 1$ 
\end{enumerate}
where $\epsilon_j\in \mathbb F_2$.

\noindent {\bf Case (1).} Assume that $\llbracket \theta_{j_0} \rrbracket_I, \ldots ,\llbracket \theta_{j_p} \rrbracket_I$ correspond to $\epsilon_{j_0}= \ldots =\epsilon_{j_p}=1$ with $j_0, \ldots,j_p\in \{ 0,\ldots, m\}$, $j_0<\ldots <j_p$, and $\epsilon_j=0$ for $j\not= {j_0}, \ldots {j_p}$. Then $C_q(\llbracket \theta_{j_0} \rrbracket_I\oplus \ldots \oplus \llbracket \theta_{j_p} \rrbracket_I)$  requires at least $q$ atoms below  $\llbracket \theta_{j_0} \rrbracket_I \oplus \ldots \oplus \llbracket \theta_{j_p} \rrbracket_I$, and this is equivalent (in every Boolean algebra) to the finite disjunction 
$$ \bigvee  (C_{\lambda_0}(\llbracket \theta_{j_0} \rrbracket_I) \wedge \ldots \wedge C_{\lambda_p}(\llbracket \theta_{j_p} \rrbracket_I))
$$
as $(\lambda_0, \ldots, \lambda_p)\in \mathbb N^{p+1}$, $\lambda_l \leq q$ for $l=0, \ldots , p$, and $\lambda_0+ \ldots + \lambda_p=q$. 

\medskip

\noindent 
{\bf Note.} 
{\rm  There can be some redundancy here, but we do not care and we do not eliminate it. The $C_{\lambda_l}$'s may be different from the original $C_d$. We stress that only finitely many conditions are required.}

\medskip

\noindent {\bf Case (2).} a) If $f_h=1$ then $C_d(f_h)$ holds if and only if $d\leq |I|$. 

b) Assume now $f_h=d_h+1$ for some polynomial $d_h$ as in Case (1) then $$f_h(\llbracket \theta_0 \rrbracket_I, \ldots ,\llbracket \theta_m \rrbracket_I)=\llbracket \theta_{j_0} \rrbracket_I \oplus \ldots \oplus \llbracket \theta_{j_p} \rrbracket_I+1=I-(\llbracket \theta_{j_0} \rrbracket_I \oplus \ldots \oplus \llbracket \theta_{j_p} \rrbracket_I)=$$ 

$$ =\bigoplus_{\substack{j=0 \\ j\not= j_0, \ldots , j_p}}^m\llbracket \theta_{j} \rrbracket_I.$$
The last equality is true since $\langle  \theta_0, \ldots, \theta_m\rangle$ is a partition sequence. So, 
$$C_r(\llbracket \theta_{j_0} \rrbracket_I \oplus \ldots \oplus \llbracket \theta_{j_p} \rrbracket_I+1) \mbox{ is equivalent to } C_r(\bigoplus_{\substack{j=0 \\ j\not= j_0, \ldots , j_p}}^m\llbracket \theta_{j} \rrbracket_I),$$ and we are in Case 1.




\

We have shown that $\Psi$ in Theorem \ref{FVsentences} is a Boolean combination of finitely many conditions $C_r(\llbracket \theta_{j} \rrbracket_I)$, where $r$ and $\theta_j $ depend only on $\phi$, and they are independent of $I$ and $\mathcal M_i$'s.

 The $\mathcal L_{B}^C$-formula $\Psi(x_0, \ldots, x_m)$  is a Boolean combination of atomic formulas in $\mathcal L_{B}^C$ uniformly for all infinite atomic Boolean algebras, and for {\it sufficiently large} finite Boolean algebras. It contains finitely many identities in the language of Boolean algebras, and finitely many conditions  $C_j(\llbracket \theta_h \rrbracket)$ and $\neg C_k(\llbracket \theta_r \rrbracket)$, saying that $\llbracket \theta_h \rrbracket$ has at least $j$ atoms below it and $ \llbracket \theta_r \rrbracket$ does not have $k$ atoms below it. In other words, $\llbracket \theta_h \rrbracket$ has at least $j$ elements,  and $\llbracket \theta_r \rrbracket$ has cardinality strictly less than $k$, respectively.  Clearly, it cannot happen that both $C_n(\llbracket \theta_j \rrbracket)$ and $\neg C_k(\llbracket \theta_j \rrbracket)$ for $k<n$ appear in $\Psi$. If $k\geq n$ then the two conditions are equivalent to $$C_n(\llbracket \theta_j \rrbracket)\vee C_{n+1}(\llbracket \theta_j \rrbracket)\vee \ldots \vee C_{k-1}(\llbracket \theta_j \rrbracket).$$ 
 
 \medskip
 
 \noindent
 {\it Proof of Theorem 2.1.  } Assume $ \prod_{i\in I}\mathcal M_i \models \phi $. Theorem \ref{FVsentences} implies that there is an acceptable sequence $\zeta =\langle  \Psi, \theta_0, \ldots, \theta_m\rangle$ such that 
 \begin{equation}
 \label{soddisfacibilità}
 \mathcal P(I) \models \Psi (\llbracket \theta_0 \rrbracket_I, \ldots ,\llbracket \theta_m \rrbracket_I). 
 \end{equation}
 

For each $i\in I$ there is a unique $j\in \{ 1, \ldots , m\}$ such that $i\in \llbracket \theta_j \rrbracket$. If $ \llbracket \theta_j \rrbracket=\emptyset$ then $\mathcal M_i\models \neg \theta_j$, for all $i\in I$.  
 Fix $j \in \{0, \ldots ,m\} $ such that $ \llbracket \theta_j \rrbracket\not=\emptyset$. The $\mathcal M_i$'s are pseudofinite structures, hence  for each $i\in \llbracket \theta_j \rrbracket$ there is a finite $\mathcal L$-structure  $\mathcal M_i^f$ such that that  $\mathcal M_i^f\models \theta_j$.

Now we consider the finite $\mathcal L$-structures $(\mathcal M_i^f)_{i\in I}$. The partition of $I$ determined by $\theta_0, \ldots, \theta_m$ and associated to $(\mathcal M_i^f)_{i\in I}$ coincides with that of $(\mathcal M_i)_{i\in I}$. So, (\ref{soddisfacibilità}) implies 
 $ \prod_{i\in I}\mathcal M_i^f \models \phi $. 
 
 Now, if $I$ is a finite set then $ \prod_{i\in I}\mathcal M_i$ is pseudofinite.  So, assume $I$ is infinite.

 In \cite{FV} the authors show that if a sentence is true in a product of structures then there is a finite product where the sentence is true. 
 This is Theorem 6.6 in  \cite{FV}, and it is proved   using ideas of Skolem (going back to 1919) which we consider a primordial elimination of quantifiers for atomic Boolean algebras. We now describe Skolem's result. 
 
 For any formula $\Phi(X_0,\ldots, X_{m})$ in $\mathcal L_B$ we can effectively associate a natural number $M$, functions $g_k: \{0, \ldots ,m\}  \rightarrow \omega$ (for $k=0, \ldots , M-1)$ and subsets $s_k$ of  $\{0, \ldots ,m\} $ (for $k=0, \ldots , M-1)$ such that if $Z_0, \ldots ,Z_{m}$ is a partition of $I$, then 
 \begin{center}
 $\mathcal P(I)\models \Phi(Z_0,\ldots, Z_{m}) $ iff  
  there is $k<M$ such that for each $j\leq m$, 
  \end{center}
 $$|Z_j|=g_k(j) \mbox{  if } j\in s_k \mbox{  and } |Z_j|\geq g_k(j) \mbox{  if } j\not\in s_k.$$
 This says that in the Boolean algebra  $\mathcal P(I)$ satisfiability of a formula by a partition of $I$ is determined by finitely many conditions on the cardinalities of the partition sets, and these cardinalities conditions are not arbitrary, there are only finitely many possibilities and  they depend on the formula $\Phi$.
 
The following result is  Theorem 6.6 of \cite{FV}.
 
 \begin{theorem*}
For any given sentence $\varphi$ we can find effectively $N\in \mathbb N$ such that if $\prod_{i\in I}\mathcal M_i\models \varphi$, then there is a set $I'\subseteq I$ having at most $N$ elements such that $\prod_{i\in I''}\mathcal M_i\models \varphi$, provided $I'\subseteq I''\subseteq I$.
 \end{theorem*}
 
 \noindent
 {\bf Sketch of the Proof:} If $\prod_{i\in I}\mathcal M_i\models \varphi$ then by Theorem \ref{FVsentences},  
 \begin{equation}
 \mathcal P(I) \models \Psi (\llbracket \theta_0 \rrbracket_I, \ldots ,\llbracket \theta_m \rrbracket_I)
 \end{equation}
 for some $\langle  \Psi, \theta_0, \ldots, \theta_m\rangle$ acceptable sequence.  
 
 By the Skolem argument   (2) is equivalent to the following statement 
 
 there are $M\in \mathbb N$, functions $g_k$ and subsets $s_k$ of $\{0, \ldots ,m\} $ with the property  that there is $k<M$ such that for each $j\leq m$, 
 
 $\bullet$ $\mathcal M_i\models \theta_j$ for exactly $g_k(j)$ elements $i\in I$ if $j\in s_k$, and 
 
 $\bullet$   $\mathcal M_i\models \theta_j$ for at least $g_k(j)$ elements $i\in I$ if $j\not\in s_k$. 
 
 In other words,  
 $$| \llbracket \theta_j \rrbracket_I |=g_k(j) \mbox{  if } j\in s_k \mbox{  and } |\llbracket \theta_0 \rrbracket_I|\geq g_k(j) \mbox{  if } j\not\in s_k.$$

Let  $n_k = g_k(0)+\ldots +g_k(m)$ for each $k\in \{ 0,\ldots M-1\}$. So, $n_k$ is the minimum number of elements of a subset of the index set $I$ which guarantees the right cardinalities of the partition sets.  Let $N=\max \{0, n_0,\ldots ,n_{M-1}\}$. For each $j\leq m$ choose a subset $Z_j\subseteq \llbracket \theta_j \rrbracket_I$ with $|Z_j| =g_k(j)$, and let $I'=Z_0\cup \ldots \cup Z_m$. Clearly, $|I'|=n_k$, and $n_k\leq N$. So,  $N$ gives an upper bound for all the possible cardinalities that may occur according to the different $k<M$. 

Take any other set $I''$ such that $I'\subseteq I''\subseteq I$. If $|\llbracket \theta_j \rrbracket_I|= g_k(j)$ then $Z_j= \llbracket \theta_j \rrbracket_I$. Hence,  $ \llbracket \theta_j \rrbracket_I \subseteq I'$, and so  $\llbracket \theta_j \rrbracket_I \subseteq I''$, and the elements in $ \llbracket \theta_j \rrbracket_I $ are the only $i$'s in $I$ such that $\mathcal M_i\models \theta_j$. So there is no risk that we may have added elements to $I''$ which change the cardinality of $\llbracket \theta_j \rrbracket_{I''}$.

If $|\llbracket \theta_j \rrbracket_I|\geq g_k(j)$ then $Z_j$ may be a proper subset of $\llbracket \theta_j \rrbracket_I$, but still it satisfies the condition of having cardinality at least $\geq g_k(j) $.  Hence,  there are at least $\geq g_k(j) $ elements $i$  in $I''$ such that $\mathcal M_i\models \theta_j$. Then $\prod_{i\in I''}\mathcal M_i\models \varphi$. 
 
 


\smallskip

We now apply Theorem 6.6 in \cite{FV}  to the family of finite structures $(\mathcal M_i^f)_{i\in I}$, and to the sentence $\phi$, and we find  a finite subset $J$ of $I$ such that  $ \prod_{i\in J}\mathcal M_i^f \models \phi $ (it is enough to choose $J=I'$). Clearly, $ \prod_{i\in J}\mathcal M_i^f $ is a finite structure, and this completes the proof that $ \prod_{i\in I}\mathcal M_i$ is pseudofinite, and the proof of Theorem \ref{productpseudo}.

\hfill $\Box$

\begin{osserva}
{\rm 
Notice that pseudofiniteness is not preserved in general under reduced product  as we show now.  Let $I$ be an infinite index set, and consider the Frech\'et filter $F$ of cofinite subsets of $I$. Consider any family of Boolean algebras $(B_i)_{i\in I}$, and let $\prod_F B_i$ be the reduced product of $(B_i)_{i\in I}$.  This is an atomless Boolean algebra, and no finite Boolean algebra is atomless.  (See also page 409 of  \cite{CK}, and \cite{frayne} where the authors study the truth in a reduced product in terms of truth in the structures.)

}
\end{osserva}

\bibliographystyle{plain}
\bibliography{ProductPseudoFinite_arxiv}
\end{document}